\documentclass[11pt]{amsart}
\newtheorem{theorem}{Theorem}[section]

\newtheorem{lemma}[theorem]{Lemma}

\theoremstyle{definition}

\theoremstyle{remark}

\numberwithin{equation}{section}

\newcommand{\CC}{\mathcal C}
\newcommand{\B}{\mathbb B}
\newcommand{\C}{\mathbb C}
\newcommand{\R}{\mathbb R}

\begin{document}

\title{ Diffeomorphisms of Stein structures}

\author{Klas Diederich  and Alexandre Sukhov}

\address{\begin{tabular}{lll}
Klas Diederich & & Alexandre Sukhov\\
Department of Mathematics & & U.S.T.L. \\
University of Wuppertal & & Cit\'e Scientifique \\
D-42097 Wuppertal & & 59655 Villeneuve d'Ascq Cedex\\
GERMANY & & FRANCE\\
 & & \\
{\rm klas@math.uni-wuppertal.de} & &
{\rm sukhov@math.univ-lille1.fr}
\end{tabular}
}

\subjclass[2000]{32H02, 53C15 \smallskip\\ \hspace*{3mm} To appear in Journal Geom. Anal.}


\begin{abstract}
We prove that a pseudoholomorphic diffeomorphism between two almost complex manifolds
with boundaries satisfying some pseudoconvexity type conditions cannot map a
pseudoholomorphic disc in the boundary to a single point. This can be viewed as an almost complex analogue
of a well known theorem of J.E.Fornaess.

\end{abstract}

\maketitle

\section*{ Main result}

A symplectic manifold admits compatible almost complex structures. Biholomorphic
invariants of these structures (such as suitable classes of pseudoholomorphic curves) allowed  Gromov
\cite{Gr} to obtain important results on symplectic rigidity.
Since his work  the study of almost complex manifolds is
rapidly increasing. Smoothly bounded domains in almost complex manifolds
play a substantial role in Gromov's approach to symplectic and contact geometry. These domains
usually arise with some additional convexity type properties which may be viewed as  almost complex analogues of the notion
of pseudoconvexity in the classical complex analysis. In particular, Eliashberg-Gromov \cite{ElGr}
and more recently, A.-L.Biolley \cite{Bi} introduced and studied under various assumptions
so called {\it Stein structures} and Stein domains
in almost complex manifolds. Since an almost complex manifold $(M,J)$ of general type does
not admit (even locally) non-constant holomorphic functions,  a Stein structure can be defined in
terms of a properly exhaustive strictly  plurisubharmonic function $\rho$ (when a symplectic form $\omega$ is fixed, a relation of type
$dd^{c}_J\rho = \omega$ is required). Sublevel sets of such a function $\rho$ are usually called Stein domains in $(M,J,\rho)$.

In this paper we prove the following

\begin{theorem}\label{wr}
Let $D$ and $D'$ be bounded domains with $C^2$ boundary in almost complex manifolds of (complex)
dimension $n$, and  assume that the target domain $D'$ admits a defining
function on a neighborhood $ U$ of the boundary $\partial D$ which is plurisubharmonic on $ U \cap D$.
Let also $F: D \longrightarrow D'$ be a pseudoholomorphic diffeomorphism
 with a $C^2$-extension $F: \overline D \longrightarrow \overline D'$. Then $F$ is a diffeomorphism between
$\overline D$ and $\overline D'$.
\end{theorem}

In the case where  almost complex structures are integrable, this is a well known result of Fornaess \cite{Fo1}.
 We state it here under a stronger assumption of plurisubharmonicity of the defining function since
the existence problem of bounded  plurisubharmonic exhaustion function similar to the Diederich-Fornaess result \cite{DiFo1}
remains  open in the almost complex category. We also point out that in the case of integrable structures
various results of this type under different assumptions have been obtained by Diederich-Fornaess \cite{DiFo2}, Fornaess \cite{Fo2},
Pinchuk \cite{Pi} in connection with Fefferman's mapping theorem. In the general almost complex case, the analogue of this theorem is known only
for strictly pseudoconvex domains  \cite{CGS,GaSu2}.

In conclusion we also add that the smoothness restrictions on the boundaries
and on the map in the hypothesis of the theorem cannot be weakened even in the integrable case (see \cite{Fo1}).

This work has been initiated when the second author visited the Max Plank Institute during the SCV colloque in September 2003.
He thanks this institution for hospitality. The work was continued while the first named author was a guest at Forschunginstitut fur
Mathematik of the ETH Zurich. He thanks this institution for giving him excellent conditions for working on it.

\section{Preliminaries}
For the convenience of the readers we briefly recall here some known facts
on the analysis and geometry of almost complex manifolds.

An almost complex structure on a smooth ($C^\infty$) real
$(2n)$-dimensional manifold $M$ is a
$C^\infty$-field $J$ of complex linear structures on the tangent
bundle $TM$ of $ M$.
We call the pair $(M,J)$ an {\it almost complex manifold}.
We denote by $J^*:T^*M \longrightarrow T^*M$ the dual of $J$ and by  $J_{st}$ the standard structure in $\R^{2n}$;
$\B$ denotes the (euclidean)  unit ball in $\R^{2n}$.

 The following frequently used lemma gives a local description of almost complex structures.

\begin{lemma}
\label{suplem1}
Let $(M,J)$ be an almost complex manifold. Then for every point $p \in
M$ and every $\lambda_0 > 0$ there exist a neighborhood $U$ of $p$ and a
coordinate diffeomorphism $z: U \rightarrow \B$ such that
$z(p) = 0$, $dz(p) \circ J(p) \circ dz^{-1}(0) = J_{st}$  and the
direct image $\hat J = z_*(J)$ satisfies $\vert\vert \hat J - J_{st}
\vert\vert_{C^3(\overline {\B})} \leq \lambda_0$.
\end{lemma}
\proof There exists a diffeomorphism $z$ from a neighborhood $U'$ of
$p \in M$ onto $\mathbb B$ satisfying $z(p) = 0$ and $dz(p) \circ J(p)
\circ dz^{-1}(0) = J_{st}$. For $\lambda > 0$ consider the dilation
$d_{\lambda}: t \mapsto \lambda^{-1}t$ in $\C^n$ and the composition
$z_{\lambda} = d_{\lambda} \circ z$. Then $\lim_{\lambda \rightarrow
0} \vert\vert (z_{\lambda})_{*}(J) - J_{st} \vert\vert_{C^3(\overline
{\mathbb B})} = 0$. Setting $U = z^{-1}_{\lambda}(\mathbb B)$ for
$\lambda > 0$ small enough, we obtain the desired statement. \qed

Denote by $T_{\C}(M)$ the complexification of the  tangent bundle $TM$ of $M$.
Then $T_{\C}M = T^{(1,0)}M \oplus T^{(0,1)}M$ where
\begin{eqnarray*}
T^{(1,0)}M = \{ X \in T_{\C}M: JX = iX \} = \{ Y - iJY, Y \in TM \}
\end{eqnarray*}
and
\begin{eqnarray*}
T^{(0,1)}M = \{ X \in T_{\C}M: JX = -i X \} = \{ Y + iJY, Y \in TM \}
\end{eqnarray*}
Using the canonical identification of the complexification $T^*_{\C}M$ of the cotangent bundle
$T^*M$ with $Hom(T_{\C}M,\C)$ we define the space $T^*_{(1,0)}M $ of complex (1,0) forms as
the space of forms $w \in T^*_{\C}M$ satisfying $w(X) = 0$ for any $X \in T^{(0,1)}M$ and
the space $T^*_{(0,1)}M $ of complex (0,1) forms as
the space of forms $w \in T^*_{\C}M$ satisfying $w(X) = 0$ for any $X \in T^{(1,0)}M$.
Every complex 1-form $w$ on $M$ may be uniquely
decomposed as $w=w_{(1,0)} + w_{(0,1)}$, where $w_{(1,0)} \in
T^*_{(1,0)}M$ and $w_{(0,1)} \in T^*_{(0,1)}M$, with respect to the structure
$J$.
This allows to define the operators $\partial_J$ and
$\bar{\partial}_J$  on
$M$~: for a complex smooth function $u$ on $M$, we set $\partial_J u =
du_{(1,0)}$ and $\bar{\partial}_Ju = du_{(0,1)}$.

Let $\Gamma = \{ r = 0 \}$ be a real smooth hypersurface in $M$.
We denote by $H^J(\Gamma)$ the $J$-holomorphic tangent bundle
$T\Gamma \cap JT\Gamma$. If  $p$ is a point of  $\Gamma$
then
$$
H_p^J(\Gamma) = \{v \in T_pM : dr(p)(v) =
dr(p)(J(p)v) = 0\} = \{v \in T_pM :
\partial_Jr(p) (v-iJ(p)v) = 0\}.
$$

The {\it Levi form} of $\Gamma$ at $p$ is the map defined on
$H^J_p(\Gamma)$  by $\mathcal L^J(\Gamma)(X_p) = J^\star dr[X,JX]_p$,
where the vector field $X$ is any section of   $H^J \Gamma$ such that $X(p) = X_p$. The hypersurface $\Gamma$  is called
 (strictly) $J$-pseudoconvex if its Levi form $\mathcal L^J(\Gamma)$
is (strictly) positive  on $H^J(\Gamma)$.  If $\rho$ is any $\CC^2$ function on $M$, the Levi
form of $\rho$ is defined on $TM$ by $\mathcal L^J(\rho)(X):=-d(J^\star d\rho)(X,JX)$.
 A $\CC^2$ real valued function $r$ on $M$ is
$J$-plurisubharmonic on $M$ (resp. strictly $J$-plurisubharmonic)
if and only if $\mathcal L^J(r)(X) \geq 0$ for every $X \in TM$ (resp.
$\mathcal L^J(r)(X) > 0$ for every $X \in TM \backslash \{0\}$).

A smooth map $f:(\tilde M, \tilde J) \longrightarrow (M,J)$ is called
pseudoholomorphic or $(\tilde J, J)$-holomorphic if its differential satisfies the following
holomorphy condition~: $df \circ \tilde J =  J \circ df$ on $TM$. In the case
$(\tilde M, \tilde J) = (\Delta, J_{st})$, where $\Delta$ denotes the unit disc of $\C$, the map $f$ is called a $J$-holomorphic
disc. An upper semicontinuous function on $M$ is called plurisubharmonic if its composition with
any $J$-holomorphic disc is subharmonic on  $\Delta$.
For $C^2$ functions this definition is equivalent to the characterization  in terms of the
Levi form giving above  (see  \cite{IvRo}).
In view of the fundamental theorem of Nijenhuis-Woolf  \cite{NiWo} on the local existence of $J$-holomorphic discs
(see Sikorav \cite{Si} for a simple proof based on the isotropic
dilations of coordinates and the implicit function theorem), any
almost complex manifold admits locally non-constant plurisubharmonic functions.
Furthermore, the notion of the Kobayashi metric can be introduced as in the integrable case (see \cite{GaSu1}, \cite{IvRo}).

\section{Proof of the theorem}
 The statement is local, so we may fix boundary points  $p \in \partial D$ and $p' \in \partial D'$ with $F(p) = p'$ and
we may choose local complex coordinates $z$ and $w$
 centered at these points satisfying the condition of  lemma \ref{suplem1}.
This reduces  the general situation to the case where $U$ and $U'$ are open neighborhoods of the origin in $\C^n$ with
standard complex coordinates $z =(z_1,...,z_n)$ and $w = (w_1,...,w_n)$ respectively, $J$ and $J'$ are  smooth real  $(2n \times 2n)$-matrix
valued functions defined on $U$ and $U'$ respectively, with $J(0) = J'(0) = J_{st}$ and $F: D \cap U \longrightarrow D' \cap U'$ is a $(J,J')$-holomorphic
diffeomorphism with a $C^2$ extension to $\partial D \cap U$. We may also assume that $D \cap U = \{ z \in U: \phi(z) < 0 \}$ and $D' \cap U' =
\{ w \in U': \rho(w) < 0 \}$. Here $\phi$ and $\rho$ are $C^2$ real functions with non-vanishing gradients on $U$ and $U'$ respectively;
furthermore, the function $\rho$ is $J'$-plurisubharmonic on $D' \cap U'$. Finally,  we may assume that
$$\phi(z) = Re z_n + O(\vert z \vert^2), \rho(w) = Re w_n + O(\vert w \vert^2).$$

Furthermore, shrinking $U$ if necessary, we may assume after  isotropic dilations ( as in the proof of lemma \ref{suplem1} )
that a map $f:\Delta \longrightarrow U$
is a $J$-holomorphic disc if and only if it satisfies

\begin{eqnarray}
\label{CR}
f_{\overline\zeta} + Q(f)\overline{(f_\zeta)} = 0
\end{eqnarray}
Here $Q$ is a smooth $(n \times n)$-matrix function $Q$ representing the matrix of the real endomorphism
$(J_{st} + J)^{-1}(J - J_{st})$ (it is easy to see that this endomorphism is antilinear with respect to the standard complex structure).
We may assume that $\parallel Q \parallel_{C^\alpha(U)} << 1$ for an arbitrary fixed {\it non-integral}
$\alpha > 2$. It is important for our proof to notice that the system (\ref{CR}) is elliptic.
A similar elliptic PDE system  can be written for $J'$-holomorphic discs with values in $U'$.

{\it Step 1: The Hopf lemma.} Fix the vector $e_n =(0,...,0,1)$ in $\C^n$. Since $J(0) = J_{st}$,
according to Nijenhuis-Woolf \cite{NiWo}  there exists a $J$-holomorphic disc $f:\Delta \longrightarrow U$
that is a solution of (\ref{CR}) of class $C^\alpha(\overline\Delta)$ satisfying $f(\zeta) = \zeta e_n + O(\vert \zeta \vert^2)$. Then the
subdomain $\Omega := f^{-1}(D \cap U) \subset \Delta$ has a $C^2$-boundary near the origin and we may apply the Hopf lemma to the
function $\rho \circ F \circ f$ negative and subharmonic on $\Omega$. This implies that $\partial F_n/\partial z_n(0) >0$.
Without loss of generality we assume  that
\begin{eqnarray}
\label{Hopf}
\frac{\partial F_n}{\partial z_n}(0) =1.
\end{eqnarray}

Since $F$ is $(J,J')$-holomorphic on $D$ and the structures $J$ and $J'$ coincide at the origin with the standard structure,
the tangent map $F'(0)$ is linear with respect to $J_{st}$. In particular,

\begin{eqnarray}
\label{tangent}
\frac{\partial F_n}{\partial z_j}(0) = 0, j=1,...,n-1
\end{eqnarray}
in view of the condition $F(\partial D) \subset \partial D'$.

{\it Step 2: Local coordinates and  discs in the target domain. }
Assume by contradiction that $F'(0)$ is degenerate. Since the holomorphic tangent space $H_0^J(\partial D)$  coincides with $\C^{n-1} \times
\{ 0 \}$,  there exist two vectors $a^0 = (a_1,...,a_{n-1},0)$ and $b^0 = (b_1,...,b_{n-1},0)$ in $H_0^J(\partial D)$ such that $a^0 \neq b^0$ but $F'(0)(a^0)
= F'(0)(b^0)$. Then the  vectors
$a = (a_1,...,a_{n-1},1)$ and $b = (b_1,...,b_{n-1},1)$ are transversal to
$\partial D$ at the origin, $a \neq b$ and in view of (\ref{Hopf}), (\ref{tangent})  $F'(0) (a) = F'(0) (b) = v = (v_1,...,v_n) \in T_0(\C^n)$ with $v _n = 1 $.
In particular, the vector $v$ is also transversal to $\partial D'$.  Consider a $\C$-linear transformation $L$ of $\C^n$ such that $L(v) = e_n$ and
$L(T_0(\partial D')) = T_0(\partial D')$. Of course, in general $L$ is not a biholomorphism for the structure $J'$, but
pushing the structure $J'$ forward by $L$, we obtain a new structure $L_*(J'): = L \circ J' \circ L^{-1}$
still equal to $J_{st}$ at the origin and such that $L$ is $(J',L_*(J'))$-holomorphic.
So we may denote again this structure  by $J'$ and suppose that $v = e_n$.

Set $\hat c = (c_2,...,c_n)$ and $e_1 = (1,0,...,0)$. Consider  a foliation of a neighborhood of the origin by a family of embedded
$J'$-holomorphic discs $f(\bullet,\hat c):\Delta \longrightarrow U'$ smoothly depending on a parameter $\hat c$ from a neighborhood of the origin in $\C^{n-1}$ and
satisfying $f(0,\hat c) = \hat c$, $f_{\zeta}(0,0) = \zeta e_1 + O(\vert \zeta \vert^2)$ (the existence of such a foliation follows from the Nijenhuis-Woolf
theorem on the existence and smooth dependence on the initial data of solutions to the equations  (\ref{CR})). Since the norm $\parallel Q \parallel_{C^\alpha}$
is small enough, we may choose the disc $f(\bullet,0)$ close enough to the disc $\zeta \mapsto \zeta e_1$ in the $C^\alpha$ norm. So for parameters $\hat c$ close enough
to the origin all the discs are embedded.  By the implicit function theorem the system of equation $w_j = f_j(\zeta,\hat c)$ , $j=1,...,n$ defines a diffeomorphism
$w \mapsto (\zeta,\hat c)$ and pushing the structure $J'$ by this diffeomorphism, we consider $(\zeta,\hat c)$ as new coordinates.
We again denote these new coordinates by $w =
(w_1,\hat w)$ and keep the notation $J'$ for the matrix representation of our almost complex structure in these coordinates. In the new system of coordinates the maps $\zeta \mapsto \zeta e_1$ are $J'$-holomorphic. If $J'(0)$ does not coincide with the standard structure, we again push $J'$ forward by a suitable linear map (and still keep the notation $J'$).  After this change of coordinates our family of discs takes the form $\zeta \mapsto l(\zeta) X + \tilde c$,
where $X$ is a vector in $T_0(\C^n)$, $l$ is a complex valued $\R$-linear function on $\C$  and the parameter $\tilde c$ runs over a neighborhood of the origin in a
complex hyperplane
of $\C^n$ transversal to $X$. However, the map $\zeta \mapsto l(\zeta) X$ is $(J_{st}, J')$-holomorphic and $J'(0) = J_{st}$. This implies that its linear part at the origin is $\C$-linear. So also the function $l$ is $\C$-linear. Finally, after a suitable $\C$-linear transformation we obtain that $l(\zeta) = \zeta$, $X = e_1$ and we still have $J'(0) = J_{st}$.
Moreover,we also may achieve by this transformation the equalities $v = e_n$ and $\rho(w) = Re w_n + O(\vert w \vert^2)$.
Thus, in the new coordinates  all our previous assumptions hold and additionally, the discs $d(\zeta,\hat w) := \zeta e_1 + \hat w$ are $J'$-holomorphic.
This implies that in these coordinates

$$
J' = \left(
\begin{array}{cll}
J_{st}^{(2)} & & C\\
0 & & B
\end{array}
\right).
$$
where $J_{st}^{(2)}$ denotes the standard complex structure of $\R^2$.

{\it Step 3: the Kobayashi metric.} Consider the rays $z = -t a$ and $z = -t b$, $t < 0$
 in $D \cap U$ for $t > 0 $ small enough. Their images $\Gamma_1(t) = F(-t a)$ and $\Gamma_2 = F(-t b)$
are  $C^2$ curves  in $\overline D'$ with the common tangent vector $e_n$ at the origin.  Since $D' \cap U' =
\{ Re w_n + 0(\vert w \vert^2) < 0 \}$ and $\partial D'$ is of class $C^2$, there exists a constant $R > 0$ such that the domain
$G: = \{ w \in U' : Re w_n + R \vert w \vert^2 < 0 \}$ is  contained in $D' \cap U'$. Clearly, $\Gamma_j(t) \in G$ for $t$ small enough.

Consider now the above family of $J'$-holomorphic discs $d_s(\bullet):= d(\bullet,\hat w)$ with $\hat w(s) = (0,...,0,-s) \in \C^{n-1}$, for $s > 0$ small enough.
Since $d_0 (\zeta) = \zeta e_1$ and $\Gamma_j$ are tangent to $e_n$ at the origin, the image $\Delta \times \{ \hat w(s) \}$ of the disc $d_s$
intersects every curve $\Gamma_j$
precisely at one point $q_j(s)$. The map $F$ is a diffeomorphism, so for $j=1,2$ there exists a unique $t_j(s) > 0$ such that
$q_j(s)= \Gamma_j(t_j(s))$ if $s$ is small enough. The Taylor expansion gives $F(ta) = - t e_n + O(\vert t \vert^2)$. So
$t_j(s) + O(\vert t_j(s) \vert^2) = s$ and  $t_j(s) \sim s$ for $j=1,2$.

Then
\begin{eqnarray}
\label{est1}
\vert q_1(s) - q_2(s) \vert \leq C_1 s^2
\end{eqnarray}
 for some positive constant $C_1$ and $ t_j(s) \sim \vert q_j(s) \vert \sim s$.
In particular,
\begin{eqnarray}
\label{est2}
\vert F^{-1}(q_1(s)) - F^{-1}(q_2(s)) \vert = \vert t_1(s) a - t_2(s) b \vert \sim s
\end{eqnarray}

 Recall also that the Kobayashi-Royden infinitesimal metric on an almost complex manifold $(M,J)$ can be defined  similarly to the case of an integrable structure.
For $p \in M$ and a tangent vector $X \in T_pM$, the value of the Kobayashi-Royden metric $K_{(M,J)}(p,X)$ is the infimum of the set of positive $\tau$ such that there exists a $J$-holomorphic
disc $f:\Delta \longrightarrow M$ satisfying $f(0) = p$ and $df(0)(\frac{\partial}{\partial (Re \zeta)}) = X/\alpha$. This definition is possible in view of the
Nijenhuis-Woolf
theorem on the local existence of pseudoholomorphic discs in any tangent direction. This metric satisfies the Schwarz lemma in the following sense:
if $\Phi: (M,J) \longrightarrow (M',J')$ is a $(J,J')$-holomorphic map, then $K_{(M',J')}(\Phi(p),\Phi'_p(X)) \leq K_{(M,J)}(p,X)$ for every point $p \in M$
and every vector $X \in T_pM$.

We may assume that the norm $\parallel J - J_{st} \parallel_{C^2(U)}$ is small enough
such that the function $\vert z \vert^2$ is strictly $J$-plurisubharmonic on $U$. Then it follows from proposition 3 of \cite{GaSu1} (see also the work of Ivashkovich-Rosay \cite{IvRo} for similar results) that  there exists  a constant $C_3 > 0$ such that
\begin{eqnarray}
\label{Kobayashi}
K_{(U,J)}(p,X) \geq C_3\vert X \vert
\end{eqnarray}
for any $p \in U$ and any tangent vector $X \in T_p(\R^{2n})$.   On the other hand, for every disc $\Delta \times \{ \hat w(s) \}$ considered above its intersection
with  $G$ is a disc of  radius $s^{1/2}$.   Let $\zeta_j(s) \in s^{1/2}\Delta$ be the points defined by $d_s(\zeta_j(s)) = q_j(s)$. Since the curves $\Gamma_j$ are tangent
to the real axis $ e_n$ at the origin and the discs $d_s$ are centered on this axes, we have
\begin{eqnarray}
\label{segment}
\vert \zeta_j(s) \vert \leq C_4 s^2
\end{eqnarray}
for some constant $C_4 > 0$.
Recall that the Poincar\'e metric of the disc $s^{1/2}\Delta$ at the point $\zeta$ on a tangent vector $X$ is given by

$$K(\zeta,X) = \frac{2^{1/2}s}{(s - \vert \zeta \vert^2)} \vert X \vert.$$
For any $s$ consider now the map $\Psi_s:= F^{-1} \circ d_s : s^{1/2}\Delta \longrightarrow D \cap U$ which is $(J_{st},J)$-holomorphic.It follows by the Schwarz lemma that
$$K(\zeta, X) \geq K_{(U,J)}(\Psi_s(\zeta), {\Psi'}_s(X)) \geq C_3 \vert {\Psi_s}'(X) \vert$$
in view of the estimate (\ref{Kobayashi}). This implies the existence of a positive constant $C_5$ such that

\begin{eqnarray*}
\parallel {\Psi_s}'(\zeta) \parallel \leq C_5 (s^{1/2} - \vert \zeta \vert)^{-1}
\end{eqnarray*}

Together with      (\ref{segment})   we get

\begin{eqnarray*}
\vert \Psi_s(\zeta_1(s)) - \Psi_s(\zeta_2(s)) \vert \leq \int_0^1 \parallel {\Psi_s}'(\zeta_1(s) + \tau(\zeta_2(s) - \zeta_1(s)) \parallel \vert \zeta_2(s) - \zeta_1(s) \vert
 d\tau \leq C_6 s^{3/2}
\end{eqnarray*}
for some constant $C_6 > 0$. On the other hand, $\vert \Psi_s(\zeta_1(s)) - \Psi_s(\zeta_2(s)) \vert \sim s$ in view of (\ref{est2}). This contradiction proves the theorem.

\end{document}